\documentclass[20pt]{article}

\usepackage[margin=1in]{geometry}
\usepackage[T2A]{fontenc}
\usepackage[utf8]{inputenc}

\usepackage{amsfonts}
\usepackage{amssymb}
\usepackage{amsmath}
\usepackage{amsthm}

\usepackage{dsfont}
\usepackage{amsmath, amsthm, amssymb, mathrsfs}
\usepackage[abbrev,non-sorted-cites]{amsrefs}
\usepackage[all]{xy}
\usepackage{graphicx}
\usepackage{extarrows}
\usepackage{hyperref}
\usepackage{xcolor}
\usepackage{tikz-cd}
\usepackage{miama}
\usepackage{calligra}
\usepackage[T1]{fontenc}

\begin{document}
\begin{Huge}
\centerline{\bf On Diophantine exponents of lattices}
\end{Huge}

\vskip+0.5cm
\begin{Large}
 
\centerline
{by   Nikolay Moshchevitin
\footnote{
Technische Universität Wien,

e-mail:nikolai.moshchevitin@tuwien.ac.at

ORCID: 0000-0001-8908-4826

The research is supported
 by Austrian Science Fund (FWF), Forschungsprojekt PAT1961524.
}}
\end{Large}
\vskip+1cm
\begin{small}
{\bf Abstract}:
We describe the spectrum of ordinary Diophantine exponents for d-dimensional lattices. The result reduces the problem to two-dimensional case and uses argument of metric theory.

\vskip+0.2cm
{\bf Keywords}: lattices, homogeneous minima, Diophantine exponents

\vskip+0.2cm
{\bf MSC}: 11H16, 11H46, 11J13
\end{small}

\vskip+1cm

\section{Spectrum of Diophantine exponents}
In the present paper we deal with Euclidean space 
  $\mathbb{R}^d$, whose points are regarded as column vectors
 $ {\bf x} = (x_1,...,x_n)^\top$.  
We define
 $$
 X=|{\bf x}| = \max_{1\le i \le d}
|x_i|
$$ 
to
be the sup-norm of 
${\bf x}$. Furthermore, we set
 $$
  \Pi ({\bf x}) = 
 \prod_{1\le i  \le d} |x_i|^{1/d}
 .
 $$
 When it is convenient to consider 
  $d$-th power of $
 \Pi ({\bf x}) $, we write $
 \Pi^d ({\bf x}) = 
 (\Pi ({\bf x}))^d$.

 For $d$-dimensional lattice $ \Lambda \subset \mathbb{R}^d$  we introduce
 an analogue of irrationality measure function,
 $$
  \psi_\Lambda (t) = \min_{{\bf x}\in \Lambda: 0<|{\bf x}|\le t} \,\,  \Pi({\bf x})
  $$
  and define the
 {\it Diophantine exponent} 
 $ \omega (\Lambda)$
by
 $$
 \omega (\Lambda) =
  \sup \{\gamma \in \mathbb{R}:\,\,\,\liminf_{t \to +\infty}
  t^\gamma \cdot     \psi_\Lambda (t)< +\infty\} 
 $$
 $$
 =
 \sup \{\gamma \in \mathbb{R}:\,\,\,
 \text{there exist infinitely many }\,\,\, {\bf x}\in \Lambda\,\,\,\text{such that }\,\,\,
  \Pi ({\bf x}) \le |{\bf x}|^{-\gamma}
 \}.
 $$
 By Minkowski's convex body theorem, for any $d$-dimensional lattice $\Lambda$
there exist infinitely many points    ${\bf x} \in \Lambda$ such that
 $   \Pi^d({\bf x}) \le {\rm covol}  \, \Lambda$. Consequently, 
  for every $d$-dimensional lattice 
  $\Lambda$ one has
 $$
  \omega (\Lambda)\ge 0.
  $$
On the other hand, there exist lattices $\Lambda$ for which
  \begin{equation}\label{0}
  \inf_{{\bf x}\in \Lambda\setminus \{ {\bf 0}\}} \Pi ({\bf x}) >0,
\end{equation}
and for such lattices we clearly have
$
  \omega (\Lambda)= 0.
$
  For example, this holds for any lattice 
  associated with a complete module in a totally real algebraic field of degree $d$.
  Moreover, it follows from a result by Skriganov \cite{SKRI} that  equality $
  \omega (\Lambda)= 0
$ holds for almost all lattices.

  \vskip+0.3cm 
  The problem of determining the set 
   $$
{\bf S}_d =
 \{ \omega(\Lambda ):\,\,\, \Lambda \,\,\, \text{is a  }\,\, d\text{-dimensional lattice} \}
 $$ 
  of all possible values of $\omega(\Lambda) $ has recently been studied in a series of papers by O. German \cite{GDEL,GIZV,GermanMonat,GUMN}. 
In particular he mentioned that in the case $ d=2$  one obviously has
$  {\bf S}_2 = [0,+\infty]$. Meanwhile, the  natural question 
of whether 
 the equality
$  {\bf S}_d = [0,+\infty]$ holds in higher dimensions remaind open. 
Using a special construction related to existence results for linear  forms of a given Diophantine type, 
German  proved (see \cite{GIZV}, Theorem 1)   
that for $ d\ge 3$ one has
$$
\left[ 3 -\frac{d}{(d-1)^2},+\infty\right] \subset  {\bf S}_d .
$$
 German  also   constructed  several other discrete sets of values from 
 $  {\bf S}_d$ using  Schmidt's subspace theorem (see \cite{GDEL}, Section 4)  . 
   \vskip+0.3cm 
   In the present paper we prove the following
   
    \vskip+0.3cm 
     {\bf Theorem 1.}
 {\it For any $ d\ge 3$ one has
$  {\bf S}_d = [0,+\infty]$.
} 
     \vskip+0.3cm 
 
 To the best of our  knowledge, the  methods of
 multiparametric geometry of numbers  \cite{GermanMonat, Omri}
do not yield such a result. Our proof  relies on a relatively simple argument from the  metric theory of Diophantine approximation  (for general results see, for example, the book by   Sprindzhuk \cite{S};
 some  related older existence results are due to  Jarn\'{\i}k \cite{J1,J2}) and  extends author's constructions from 
  \cite{M1,M2}. 
  
    \vskip+0.3cm 
    In the next section we formulate a more general Theorem 2 from which Theorem 1 follows and discuss possible further improvements. The rest of the paper deals with a proof of Theorem 2.

 \section{General proposition and some remarks}
 
 Theorem 1 is an immediate corollary of the following
 
   \vskip+0.3cm 
     {\bf Theorem 2.}
     {\it Let $d\ge 3$. Consider real function $w(t)$ decreasing to  zero as $ t \to +\infty$.
    Assume that $w(t)$ satisfies a technical condition    
       \begin{equation}\label{www00}
  \lim_{\delta \to 0} \sup_{t\ge 1} \frac{w(t)}{w((1+\delta)t)} =1.
      \end{equation}
     Assume that  the series 
 \begin{equation}\label{www}
 \sum_{\nu=1}^\infty \nu^{d-1} w \left({2^\nu}\right) 
 \end{equation}
 converges. Then, for every positive $\varepsilon$, there exists a lattice $\Lambda$ such that there are infinitely many points ${\bf x}\in \Lambda$ satisfying
 \begin{equation}\label{q}
  \Pi^d ({\bf x})  \le (1+\varepsilon) w (|{\bf x}|),
 \end{equation}while  the stronger inequality
 \begin{equation}\label{p}
  \Pi^d ({\bf x})  \le  (1-\varepsilon) w(|{\bf x}|)
 \end{equation}
 holds for at most finitely many points ${\bf x}\in \Lambda$.}
   \vskip+0.3cm 
    
   To get Theorem 1 from Theorem 2,
   it suffices to
    apply it to functions $w(t) = t^{-\gamma d}$ for  $\gamma >0$. 
    Moreover, note
   that for any $\eta>0$ function $\omega(t) = (\log (1+t))^{-d-\eta}$  also satisfies the assumptions
    of Theorem 2. Hence,   there exists a lattice $\Lambda = \Lambda(\eta,\varepsilon)$  such that  there are
   infinitely many points ${\bf x}\in \Lambda(\eta,\varepsilon)$ satisfying
 \begin{equation}\label{qw}
  \Pi^d ({\bf x})  \le (1+\varepsilon)   (\log (1+|{\bf x}|))^{-d-\eta},
 \end{equation}
 while  inequality 
 \begin{equation}\label{qww}
 \Pi^d ({\bf x})  \le (1-\varepsilon)   (\log (1+|{\bf x}|))^{-d-\eta}
 \end{equation}
holds only for at most  finitely many  points.
 
    \vskip+0.3cm 
    
 The main purpose of the present paper is  to obtain an optimal  qualitative result about spectrum ${\bf S}$. However, we believe that quantitative statement of Theorem 2 is not optimal.
 First, we expect that condition  related to convergence of the series (\ref{www}) can be weakened. 
 Second, we 
 are sure that a sharper existence result
 should hold, in which the constant
  $\varepsilon$  is replaces by a function  $\varepsilon (|{\bf x}|)$  decreasing  to zero
  sufficiently 
   fast as $ |{\bf x}|\to \infty$. 
   Below we briefly explain 
  how such improvements might be obtained.
   
       \vskip+0.3cm

      The assumption about the convergence of  the series (\ref{www})  is used in our proof to show that "the most part" of points ${\bf x}$ of the lattice $\Lambda$ admit a good lower bound for the value of
      $
  \Pi({\bf x})
     $. 
  The possibility of obtaining   a stronger statement 
    is related to the result by Skriganov mentioned in the previous section.  
  For completeness, we state here a precise formulation.
  Namely, we give here an exact formulation of one of the results from  \cite{SKRI},
  adapted to the notation of the present paper.
  This statement may be of importance in accordance with our topic.

  \vskip+0.3cm
{\bf Theorem A } (Lemma 4.3 from \cite{SKRI}). {\it 
 Let function $ \varphi (t): \mathbb{R}_+\to \mathbb{R}_+$ satisfy condition
  $\varphi(t) = O(\exp (at)), t\to +\infty$ with some $ a\ge 0$ and
  $$
  \int_1^{+\infty} \frac{{\rm d}t}{ \varphi (t)} < \infty.
  $$
  Let $\Lambda$ be an arbitrary $d$-dimensional lattice in $\mathbb{R}^d$.
  Then for almost all $ U \in {\rm SO}(d)$ (with respect to Haar measure on $ {\rm SO}(d)$) 
  for the lattice $U\Lambda$ 
   one has
  $$
  \liminf_{t\to +\infty} \,\,\,
   (\log t)^{d-1} \, \varphi (\log\log t) \cdot   (\psi_{U\Lambda} (t))^d>0.
  $$
  In particular, for any positive $\eta$, for almost all $U$ one has
   $$
  \liminf_{t\to +\infty} \,\,\,
   (\log t)^{1-1/d +\eta}  \cdot   (\psi_{U\Lambda} (t))>0,
  $$
  or
\begin{equation}\label{pwww}
  \Pi^d ({\bf x}) \gg_{U,\varepsilon} (\log(1+|{\bf x}|))^{1-d-\eta}\,\,\,\,
  \text{when}
  \,\,\,\,
  {\bf x} \in U\Lambda \setminus \{{\bf 0}\}.
  \end{equation}
  }
 
 \vskip+0.3cm 
We see that the exponent in the power of logarithm in (\ref{pwww}) is better that that in 
(\ref{qw},\ref{qww}).
  \vskip+0.3cm

 The second opportunity to improve Theorem 2 
 is to obtain existence results about exact order of approximation,
 that is to obtain a statement where inequalities 
 (\ref{q},\ref{p}) are replaced by precise asymptotics of the form
 $$
   \Pi^d ({\bf x})  = w  (|{\bf x}|) (1+o(|{\bf x}|)),\,\,\,\,\,|{\bf x}|\to \infty
 $$ 
 with an explicit expression in the remainder $o(|{\bf x}|)$.
 In this direction the following very precise result which was recently proven by
 Baker and Ward \cite{BW} may be useful.

\vskip+0.3cm
{\bf Theorem B} (Theorem 1.1 from \cite{BW}).  {\it
Let function $\psi(t): \mathbb{R}_+\to \mathbb{R}_+$ satisfy the condition $ t\cdot \psi(t)  \downarrow 0, t \to \infty$.
Define
$$  
\psi_1(t) = \psi(t)  (1-\omega(t)),\,\,\,\,\,
\text{where}
\,\,\,\,\,
\omega(t) = 90 t  \psi\left( \frac{1}{2\psi(t)}\right) \to 0, \,\,\,t\to \infty.
 $$
 Then there exists  
 an uncountable set $E[\psi,\omega]$ of numbers $ \alpha \in [0,1]\setminus \mathbb{Q}$  
 with the following property: for every $\alpha \in E[\psi,\omega]$,
 
 \vskip+0.3cm

 {\rm 
 ({\bf i})} there exists a sequence of continued fraction convergents  $p_{\nu_k}/q_{\nu_k}, k = 1,2,3,...$ to $\alpha$ with 
 $$
 |q_{\nu_k}\alpha - p_{\nu_k}| \le \psi (q_{\nu_k});
 $$
 
 {\rm  ({\bf ii})}  there exists $q_0$ 
 such that  
 \begin{equation}\label{BW1}
 |q\alpha - p| \ge \psi_1(q),\,\,\,\,\,
 \forall\, q \ge q_0, \,\,\forall\, p. 
\end{equation}
 In particular
 \begin{equation}\label{BW2}
|q_{\nu_k}\alpha - p_{\nu_k}|  = ||q_{\nu_k} \alpha || = \psi(q_{\nu_k})(1 + O(\omega (q_{\nu_k}))).
\end{equation}

 }

Here and everywhere in this paper we use the notation
$||\alpha|| =\min_{a\in \mathbb{Z}}|\alpha - a|$ for the distance to the nearest integer.
 It is clear that the construction from  \cite{BW}
 ensures that the set   $E[\psi,\omega]$  is dense in $[0,1]$.
  
  \vskip+0.3cm

 We should notice that  precise existence results for simultaneous approximations by Akhunzhanov \cite{a}
 may also be useful in obtaining a result
  which is more sharp than Theorem 2.

 \vskip+0.3cm

 Also we would like to discuss technical condition  (\ref{www00})
 It follows from technical condition  (\ref{www00})  that 
  \begin{equation}\label{www0000}
 0< \sup_{t\ge 1} \frac{w(kt)}{w(t)} <+\infty\,\,\,\,\,\text{for all}\,\,\,\,\, k >0.
      \end{equation}
 Technical conditions (\ref{www00})  and  (\ref{www0000}) are quite natural in statements of such a type (see, for example,  formulations of classical results by Jarn\'{\i}k  \cite{j41},
   Theorem 2 from \cite{MMJ},  or discussion in Section 2.2 from \cite{NM}). However the result of Theorem 2 should be valid under much weaker restrictions.
   
 \vskip+0.3cm
 Finally, in this section we introduce  notation
 \begin{equation}\label{pao}
\frak{P} [\varphi]  =
\left\{ {\bf x}\in \mathbb{R}^d: \,\,\,
 \max_{j=1,...,d} |x_j| \ge 1,
\,\,\,
\Pi^d ({\bf x}) \le \varphi \left( \max_{j=1,...,d} |x_j| \right)
\right\}.
\end{equation}
Our approach is related to counting points of a random lattice in this set.

 \section{ Two-dimensional subspace}
  \label{el}


\noindent
In
$\mathbb{R}^d$  with coordinates $x_1,...,x_d$ we consider  $(d-1)$-dimensional coordinate subspaces 
$$
\frak{R}_ j = \{ {\bf x} \in \mathbb{R}^d: \,\, x_j = 0\},\,\,\, j = 1,...,d,
$$
and one-dimensional coordinate subspaces 
$$
\frak{r}_i= \frak{R}_i^\perp =\bigcap_{j: \, j \neq i} \frak{R}_ j .
$$
We should note that 
$
\frak{r}_i \subset \frak{R}_ j, \,\,\, \forall i: i \neq  j.
$

\vskip+0.3cm
In our proof we fix  a two-dimensional linear subspace 
$
\frak{L} \subset \mathbb{R}^d
$
satisfying 
\begin{equation}\label{satisa}
\frak{L} \cap \frak{R}_i\cap \frak{R}_j = \{{\bf 0}\}\,\,\, \forall\, i,j = 1,...,d : i \neq j.
\end{equation}
It is clear that  a generic two-dimensional subspace satisfies  
this condition. We should note that 
\begin{equation}\label{satisb}
 \frak{r}_i\not\subset \frak{L}\,\,\, \forall\, i = 1,...,d.
\end{equation}
for every $\frak{L}$ satisfying (\ref{satisa}).
For fixed $\frak{L}$ consider one-dimensional subspaces
$$
\frak{l}_j= \frak{L}\cap \frak{R}_ j ,\,\,\,\,\, j = 1,...,d.
$$
It follows from property (\ref{satisa}) that
\begin{equation}\label{comti0}
\frak{l}_j  
\cap\frak{R}_i =\{\pmb{0}\}
 ,\,\,\,\,\, \forall \,  i\neq j ,
\end{equation}
and in particular
\begin{equation}\label{comta}
\frak{l}_j  \neq \frak{r}_i ,\,\,\,\,\, \forall \,  i, j  = 1,...,d.
\end{equation} 
We define
$
\frak{l}^\perp\subset \frak{L}
$
 as one-dimensional subspace of $\frak{L}$ which is orthogonal to 
$\frak{l}
$.
From  (\ref{satisb})  we see  that 
\begin{equation}\label{comti}
\frak{l}_1^\perp \neq \frak{r}_i ,\,\,\,\,\, \forall \,  i = 1,..,d.
\end{equation}
As $ \frak{l}_j \subset \frak{R}_j$ we see that  $ \frak{l}_j \perp \frak{r}_j$,
meanwhile
\begin{equation}\label{satisaq}
 \frak{l}_j^\perp \not\subset \frak{R}_j.
\end{equation}
By 
$\varphi_j\in (0,\pi/2)$ we denote the angle between  $ \frak{l}_j^\perp$ and $\frak{r}_j$.
Notice that all the  angles  $\varphi_j$ as well as  parameters $\pmb{\xi}, e_{i,j},  c_j $ and $c$  which we define later in this section depend on our choice of subspace $\frak{L}$.
All the constants in symbols $\ll, \gg, \asymp, O(\cdot) $ below may  also depend on $\frak{L}$.
If the constants in these symbols depend on something else, we outline this dependence.
 
\vskip+0.3cm
 Now  we take two vectors 
$
\pmb{e}_1 , \pmb{e}_2 \in \frak{L}$ of unit Euclidean length
such that
$$\frak{l}_1^\perp = \langle \pmb{e}_1 \rangle_\mathbb{R},\,\,\,
\frak{l}_1  = \langle \pmb{e}_2 \rangle_\mathbb{R}.
$$
Vectors $
\pmb{e}_1 , \pmb{e}_2 $ are defined up to sign $\pm$. For convenience we choose the sign a little bit later.
So every point  ${\bf x} \in \frak{L}$ may be represented as a sum
\begin{equation}\label{Punkt}
{\bf x} = (x_1,...,x_d)^\top = z_1\pmb{e}_1 + z_2\pmb{e}_2  \in \frak{L},
\end{equation}
We denote
$$
{\bf z}_- = (z_1,z_2)^\top, \,\,\,\,\,\,\,\, Z_- = |{\bf z}_-| = \max_{ 1\le i \le 2}|z_i|,
$$
so by triangle inequality 
\begin{equation}\label{PunktX}
X \le |z_1|+|z_2|\le 2 Z_-.
\end{equation}
In coordinates $x_1,...,x_d$ vectors $\pmb{e}_1,\pmb{e}_2$
are represented as 
\begin{equation}\label{eee}
\pmb{e}_1=
\left(
\begin{array}{c}
e_{1,1}
\cr
e_{1,2}
\cr
\vdots
\cr
e_{1,d}
\end{array}
\right),
\,\,\,
\pmb{e}_2=
\left(
\begin{array}{c}
0
\cr
e_{2,2}
\cr
\vdots
\cr
e_{2,d}
\end{array}
\right),
\,\,\,
e_{1,1}
>0
,\,\,\,\,\, e_{2,j} \neq 0,\,\,\ j = 2,...,d;
\end{equation}
here 
(\ref{satisaq})
ensures the inequality  $e_{1,1} \neq 0$ and the choice of sign for $
\pmb{e}_1$ ensures positivity of coordinate $e_{1,1}$, meanwhile
 $e_{2,1} = 0$ as
$\frak{l}_1  = \langle \pmb{e}_2 \rangle_\mathbb{R}  = \{ \pmb{x}\in \frak{L}:\,\,\, z_1 = 0\}\subset \frak{R}_1$ and non-vanishing of $e_{2,j}$  follows from  (\ref{comti0}).
Vectors $\pmb{e}_1, \pmb{e}_2$ are independent. So there exist $i\neq j$ such that 
$$
\left(
\begin{array}{c}
x_{i}\cr 
x_{j} 
\end{array}
\right)=
\left(
\begin{array}{cc}
e_{1,i} & e_{2,i}\cr
e_{1,j}& e_{2,j}
\end{array}
\right)
\left(
\begin{array}{c}
z_{1}\cr 
z_{2} 
\end{array}
\right),\,\,\,\,\,\,
\left|
\begin{array}{cc}
e_{1,i} & e_{2,i}\cr
e_{1,j}& e_{2,j}
\end{array}
\right| \neq 0
$$
and
\begin{equation}\label{dete0}
Z_- \asymp \max(|x_i|, |x_j|) \le X.
\end{equation}
Combining (\ref{PunktX}) and (\ref{dete0}) we get 
 \begin{equation}\label{dete}
Z_- \asymp X.
\end{equation}
Moreover,
 \begin{equation}\label{deteX}
X =   |z_2|\cdot \max_{2\le j \le d}|e_{2,j}|   + O(|z_1|).
\end{equation}

We define real numbers  $\xi_2,...,\xi_d, $ by the condition
$$
\frak{l}_j  = \{ {\bf x} \in \frak{L}:\,\,\, z_2 = \xi_j z_1\},\,\,\,\,\,\,
j = 2,...,d.
$$

It follows from our definitions that  for  ${\bf x} \in \frak{L}$ the  coordinates from (\ref{Punkt}) satisfy
\begin{equation}\label{X1}
|x_1| = c_1  |z_1| ,
\,\,\,\,\,\text{where}
\,\,\,\, c_1= {\cos \varphi_1 }\neq 0
\end{equation}
and
\begin{equation}\label{X2}
|x_j | =
{\rm dist}\, (\frak{l}_j,{\bf x} ) \cdot  \cos \varphi_j =
c_j |z_2 - \xi_jz_1| ,\,\,\, j = 2,...,d,
\,\,\,\,\,\text{where}
\,\,\,\, c_j = \frac{\cos \varphi_j }{\sqrt{1+\xi_j^2}} \neq 0
\end{equation}
(we use ${\rm dist}\, (\cdot,\cdot ) $ for Euclidean distance between the sets).
Now 
for the product
$\Pi ({\bf x})$ we get equality 
\begin{equation}\label{X3}
\Pi^d ({\bf x})
  = \prod_{j=1}^d |x_j| = 
c|z_1| \prod_{j=2}^d |z_2 - \xi_jz_1|\,\,\,\,\,
\text{with}
\,\,\,\,\,
c = \prod_{j=1}^d c_j.
\end{equation}

For  the  set
$$ \frak{P} = \frak{P}[\varphi_0 ]
=\{{\bf x}\in \mathbb{R}^d:\,\,\, \Pi^d ({\bf x})\le 1\}
$$ defined by the 
 function $\varphi_0 (t) =1 ,\, \forall\, t$ and for
 our subspace $\frak{L}$ we should study the properties of the intersection
$$
\frak{P}_{\frak{L}} = \frak{P} \cap  {\frak{L}}.
$$
It is clear that 
$${\bf x} \in  \frak{L} \setminus 
\frak{P}_{\frak{L}}\,\,\,\,\,
\Longrightarrow
\,\,\,\,\,
\Pi ({\bf x})>1.
$$

\vskip+0.3cm
{\bf  Lemma 1.}
{\it For any  ${\bf x} \in 
\frak{P}_{\frak{L}}$ there exists  index $j_*$ such that 
\begin{equation}\label{lema1}
|x_{j_*}| \ll \frac{1}{X^{d-1}}\,\,\,\,\,\,\,\text{and}\,\,\,\,\,\,\,\,
|x_j| \gg X\,\,\,\,  \forall\, j \neq j_*.
\end{equation}
}
\vskip+0.3cm

Proof.
Lines $\frak{l}_j, j =1,...,d$ divide two-dimensional subspace $\frak{L}$ into $2d$ connected angular domains. Let denote these domains as $D_k \subset \frak{L}, k=1,...,2d$. For each of these domains we consider the intersections
$$
D_k^0 = D_k \cap \{ {\bf x}\in  \mathbb{R}^d: \,\,\, \Pi ({\bf x}) \le 1\}.
$$
Each of these intersections is a connected domain with boundary $  \partial 
D_k^0 $
which contains of two rays $ l_{k}^+,l_k^-$ which are parts of some of lines $\frak{l}_j$ and a curve
$$
\kappa_k = D_k \cap \{ {\bf x}\in \mathbb{R}^d:\,\,\,\Pi ({\bf x}) = 1\}
$$
 which lies between $ l_{k}^+$ and $ l_k^-$, so 
 $  \partial 
D_k^0 = \kappa_k \cup  l_{k}^+\cup l_k^-$. 
Any  ${\bf x}' = (x_1',...,x_d') \in \kappa_k$ satisfies the equality $\Pi ({\bf x}') = 1$.
Consider  a point ${\bf x}' = (x_1',...,x_d') \in \kappa_k$ which continuously tends to infinity remaining in $\kappa_k$. 
Then there exists
$j_*$ such that the Euclidean distance 
${\rm dist}\, ({\bf x}', \frak{l}_{j_*})$ 
between ${\bf x}' $ and $\frak{l}_{j_*}$ tends to zero.
By (\ref{comti0})
this means that for all other coordinates $ x_j', j\neq j_*$ one has
\begin{equation}\label{due}
 |x_j'|\asymp |{\bf x}'|
 ,
 \,\,\,\,\, 
 j
 \neq j_*.
 \end{equation}
As in our consideration $\Pi({\bf x}') = 1$, we conclude that 
\begin{equation}\label{due1}
{\rm dist}\, ({\bf x}', \frak{l}_{j_*})\asymp |x_{j_*}'|
=
\frac{1}{ \prod_{j \neq j_*} |x_j'|}
\asymp \frac{1}{|{\bf x}'|^{d-1}}
\end{equation}
  as ${\bf x}'$ goes to infinity.
Now if we assume that ${\bf x} \in \frak{P}_{\frak{L}}$ and $X$ is large, there exists  $j_*, k$  and a point 
${\bf x}' \in\kappa_k $  close to ${\bf x}$ such that 
$ {\rm dist}\, ({\bf x}, \frak{l}_{j_*})\le 
{\rm dist}\, ({\bf x}', \frak{l}_{j_*})
$
and $ |x_j'| \asymp  |x_j| $ for every $ j \neq j_*$.
  This together with (\ref{due},\ref{due1}) proves Lemma 1. $\Box$

\section{Two-dimensional lattice}\label{L}
 We consider $\mathbb{R}^2$ with coordinates $z_1,z_2$.
For two real numbers 
$
\alpha,\beta \in [0,1/2]
$
we define the  lattice
$$
\Gamma =  \Gamma_{\alpha,\beta} =
\left\{
{\bf z}_{m,n} =
\left(
\begin{array}{c}
z_1 \cr
z_2\end{array}
\right)=
\left(
\begin{array}{c}
m\alpha +n \cr
m+n\beta \end{array}
\right),\,\,
m,n \in 
\,
\mathbb{Z}\right\}=
 A
\mathbb{Z}^2\,\,\,\,\,
\text{where}
\,\,\,\,\,
A=
 \left(
\begin{array}{cc}
\alpha & 1 \cr
1&\beta
\end{array}
\right).
$$
Notice that 
as 
${\rm det } \,A=
\alpha\beta -1
 $ bounded from zero,
for any 
${\bf z}= {\bf z}_{m,n}\in \Gamma$   we have 
\begin{equation}\label{mn}
 \max(|m|,|n|) \asymp |{\bf z}_{m,n}|.
 \end{equation}

\vskip+0.3cm
\noindent
Let real-valued function $ \psi (t) $ 
decreases to zero as $ t \to +\infty$.
 For an integer $n$ we write
$\psi_n = \psi(|n|)$.

\vskip+0.3cm
{\bf Lemma 2.}  {\it Let 
$
\alpha  \neq 0
$ and
real numbers $ \xi_2,...,\xi_d$ satisfy 
\begin{equation}\label{nee}
 \xi_j\alpha \neq 1,\,\,\,\,\, j = 2,...,d
.
\end{equation}
Suppose that  the 
series $\sum_{n} \psi_n$ converges.
Then   
there is a set $\mathcal{B}\subset  [0,1/2]$ of full 
 Lebesgue measure  such that for all $ \beta \in \mathcal{B}$ there exists 
 $ c = c(\beta)>0$ with 
\begin{equation}\label{rezu}
|z_2-\xi_j z_1 |> c\psi_n
,
\,\,\,\,
\forall \, {\bf z} = {\bf z}_{m,n} \in \Gamma_{\alpha, \beta} 
\,\,\,\text{with}
\,\,\,
n\neq 0\,\,\,\,\text{and all}\,\,\,\,
  j = 2,...,d.
\end{equation}
}

Proof.
Let $n\neq 0 $ and ${\bf z}\in \Gamma$. We deal with  inequalities
\begin{equation}\label{1}
 |z_2-\xi_j z_1 | <\psi_n,\,\,\,\,\,
 j = 2,...,d
\end{equation}
as with conditions on $\beta$. Consider the set 
$$
J_{m,n} =
\{\beta \in [0,1/2]:\,\, \min_{2\le j \le d} |z_2-\xi_j z_1 | <\psi_n\}=
\left\{\beta\in [0,1/2]:\,\, \exists \, j\,\,\,\text{such that} \,\,\,\left|\beta - \frac{m(\xi_j\alpha-1)}{n} -\xi_j\right| <
\frac{
\psi_n}{|n|}\right\}
$$
and the union
$$
J_{n} = \bigcup_{m: 
\frac{m(\xi_j\alpha-1)}{n} -\xi_j \in
\left[-\frac{\psi_n}{|n|}, \frac12+\frac{\psi_n}{|n|} \right]}\,\,\, J_{m,n}. 
$$
For the Lebesgue measure of $J_{m,n} $ we have a bound
$$
\mu (J_{m,n}) \le \frac{2(d-1) \psi_n}{|n|}
$$
and so  by (\ref{nee}) we conclude that 
$$
\mu (J_{n}) \ll_{\alpha} { \psi_n}
.$$

If there exists infinitely many ${\bf z}$ with (\ref{1})
then
$$
\beta \in \bigcap_{n_0=1}^\infty \bigcup_{n > n_0}  J_{n}.
$$
Now the conclusion of the lemma follows by  Borel-Cantelli argument.$\Box$

\section{ Embedding of $\Gamma$ into $\mathbb{R}^d$}\label{s5s}

We 
take arbitrary irrational $\alpha$ satisfying (\ref{nee}) and produce by Lemma 2 set  $\mathcal{B}$.
Then we take arbitrary $\beta \in \mathcal{B}$
and  consider lattice
$$
\Gamma_{\frak{L}}
=
 \left\{
 z_1\pmb{e}_1+z_2\pmb{e}_2:\,\,\,
\left(
\begin{array}{c}
z_1 \cr
z_2\end{array}
\right)=  \pmb{z}_{m,n}\in \Gamma
\right\}
\subset \frak{L}$$
 congruent to $\Gamma$ from Section \ref{L}. Now we analyse the values of $\Pi ({\bf x})$ for points 
 ${\bf x} $ from $
\Gamma_{\frak{L}}$. Recall that we use the notation $X =|{\bf x}| = \max_{j = 1,...,d} |x_j|$.

\vskip+0.3cm
{\bf  First Main Lemma.} {\it 
Consider irrational numbers $ \alpha ,\beta$, assume that condition (\ref{nee}) is satisfied and assume  
that $ \beta \in \mathcal{B}$ satisfies 
\begin{equation}\label{betalok}
 \alpha\beta \neq 1,\,\,\,\,\,\,\,\,\,\,
 \beta \neq \xi_j,\,\,\, j = 2,...,d.
 \end{equation}
Let  $p_\nu/q_\nu, \nu = 1,2,3,...$  be convergent fractions for  $\alpha$.

\vskip+0.3cm

 {\rm 
 ({\bf i})}
 Let $ {\bf x} = \pmb{z}_{m,n}\in \Gamma_{\frak{L}}$ be of the form
 $$
 (m,n) \ne s(q_\nu, - p_\nu),\,\,\, s \in\mathbb{Z}
\setminus \{0\}. 
 $$
 Then
 \begin{equation}\label{ito}
\Pi^d ({\bf x})   \gg_{\alpha,\beta }
\min
\left(1,\psi (HX)  X^{d-1}\right),
\end{equation}
with some constant $H=H(\frak{L})$ depending on $\frak{L}$.

\vskip+0.3cm

  {\rm 
 ({\bf ii})}
 Let $ {\bf x} = \pmb{z}_{m,n}\in \Gamma_{\frak{L}}$ be of the form 
 $$
 (m,n) = s(q_\nu, -p_\nu),\,\,\, s \in\mathbb{Z}.
 $$
Then
 \begin{equation}\label{ito1}
\ \Pi^d ({\bf x})
 =
 \sigma |1-\alpha\beta|^{d-1}
 |s|^d  q_\nu^{d-1}||q_\nu\alpha|| \left(1 + O\left(  \frac{||q_\nu\alpha||}{q_\nu}\right)\right) \,\,\,\,\,
\text{where}
\,\,\,\,\,
\sigma =  \sigma (\frak{L}) = e_{1,1}\prod_{j=2}^d|e_{2,j}| \neq 0,
 \end{equation}
   and $e_{i,j}$ are coordinates of vectors from (\ref{eee}).
 }
  
\vskip+0.3cm

Proof. 
Let us prove statement  ({\bf i}).
If $ x\not\in \frak{P}_{\frak{L}}$ there is nothing to prove.

If $m=0$  or $n=0$ point $\pmb{z}_{m,n}$ has coordinates
$z_1=n, z_2= n\beta$ or $ z_1 = m\alpha, z_2 = m$ respectively. We take into account conditions (\ref{betalok},\ref{nee})  and formulas
(\ref{X1},\ref{X2}) to see $|x_j|\asymp X$  for every $j$.  
So in this case we get 
$\Pi^d ({\bf x}) \gg X^d\ge 1$.
Let $mn\neq 0$, then
in the case ({\bf i})  Legendre's theorem, formula (\ref{X1}) and inequalities (\ref{mn},\ref{dete})  lead to
the inequality
\begin{equation}\label{possa}
|x_1|\asymp |z_1|= | m\alpha +n|\ge \frac{1}{2|m|}\gg\frac{1}{\max(|m|,|n|)} \asymp
\frac{1}{|\pmb{z}_{m,n}|}\asymp\frac{1}{X}
.
\end{equation}
For all other $j =2,...,d$ by formulas (\ref{X2}) and  condition (\ref{rezu}) we get
\begin{equation}\label{possa1}
|x_j|\gg_{\alpha,\beta} \psi_n.
\end{equation}
We consider index $j_*$ from  Lemma 1.
In our case for $ {\bf x} \in \frak{P}_{\frak{L}}$  equality $j_* = 1$ is not possible  for large $X$ because of
(\ref{possa}) and the first inequality from (\ref{lema1}).
If 
 $ {\bf x} \in \frak{P}_{\frak{L}}$  and  $j_* \neq  1$ we apply (\ref{possa1}) for $ j = j_*$ and 
 the second inequality from (\ref{lema1}) for all other $j\neq j_*$.
 In such a way we obtain
 $\Pi^d ({\bf x}) \gg_{\alpha,\beta}  \psi_nX^{d-1}$. As
 $ |n|\le \max(|m|,|n|) \asymp |\pmb{z}_{m,n}|\asymp X$, we 
  get (\ref{ito}).
  
  \vskip+0.3cm

Now we prove statement ({\bf ii}).  From formula (\ref{X3}) and our definitions immediately follows equality
$$
\Pi^d ({\bf x})   = 
 |s|^d \cdot
 e_{1,1} ||q_\nu \alpha||
 \cdot\prod_{j=2}^d
 |e_{1,j} (q_\nu\alpha- p_\nu) +e_{2,j}(q_\nu - \beta p_{\nu})|.
 $$
 But 
   $$
 q_\nu -\beta p_\nu = (1-\alpha\beta) q_\nu + O(||q_\nu\alpha||).
 $$
 From last two formulas asymptotic equality 
 (\ref{ito1})  follows.$\Box$

\vskip+0.3cm

At the end of this section we should note that for ${\bf x}$ of the form
${\bf x}=\pm \pmb{z}_{{q_\nu},-p{_\nu}}$ from the case  {\rm 
 ({\bf ii})} of the First Main Lemma  with $ s = \pm 1$,  by (\ref{deteX}) we have
 \begin{equation}\label{1x}
 X =   |z_2|\cdot \max_{2\le j \le d}|e_{2,j}|   + O(|z_1|)=
 q_\nu |1-\alpha\beta| \max_{2\le j \le d}|e_{2,j}| + O_\beta (||q_\nu \alpha||).
\end{equation}

\section{Multi-dimensional domains}

Let $ \varepsilon >0$ and $ X >1$.
We begin this section with covering of $(d-1)$-dimensional domain 
 $$
\frak{P}_{\varepsilon, X} 
=
\left\{ \underline{\bf x} = (x_2,...,x_d)\in \mathbb{R}^{d-1}:\,\,\,\,
|x_2\cdots x_d|\le \varepsilon,\,\,\,\ \max_{j =2,...,d} |x_j|\le X\right\}
$$
by boxes.
Let
$${\bf t} = (t_1,...,t_d) = (t_1, \underline{\bf t}) \in \mathbb{Z}^d,\,\,\,\,\,\,\,\,
\underline{\bf t} = (t_2,...,t_d)\in \mathbb{Z}^{d-1},
$$
\begin{equation}\label{suum}
t_2+...+t_d = 0.
\end{equation}
Consider parameter
$$
T_{\varepsilon, X}
=
\left[ (d-2)
\frac{\log X+ \log \left( \varepsilon^{-1/(d-1)}\right)}{\log 2} + \frac{1}{2}\right].
$$
and boxes
$$
\frak{B}_{\underline{\bf 0}} (\varepsilon) = 
\left\{\underline{\bf x} =(x_2,...,x_d)\in \mathbb{R}^{d-1}:\,\,\,
\max_{2\le j\le d} |x_j|
\le 2^{(d-1)/2}\varepsilon^{{1}/({d-1})}\right\},\,\,\,\,
\frak{B}_{\underline{\bf t}} (\varepsilon) =
G_{\underline{\bf t} } \frak{B}_{\underline{\bf 0}}(\varepsilon) 
,$$
$$\text{where}\,\,\,\,\,
G_{\underline{\bf t}}
=
\left(
\begin{array}{ccc}
2^{t_2}& \cdots& 0\cr
\cdots& \cdots&\cdots\cr
0& \cdots& 2^{t_d}
\end{array}
\right)\,\,\,\,\, \text{is} \,\,\,\,\,
(d-1)\times (d-1) \,\,\,\text{matrix}.
$$

\vskip+0.3cm

{\bf Lemma 3.}
{\it  For any  $ \varepsilon \in (0,1)$ and $ X >1$  one has 
\begin{equation}\label{geho}
\frak{P}_{\varepsilon, X} \subset \bigcup_{\underline{\bf t}}\frak{B}_{\underline{\bf t}} (\varepsilon) ,
\end{equation}
where the union is taken over all 
$\underline{\bf t}\in \mathbb{Z}^{d-1}$ such that  (\ref{suum})  holds and 
\begin{equation}\label{maax}
\max_{j=2,...,d-1} |t_j| \le T_{\varepsilon, X}
\end{equation}
}

 \vskip+0.3cm

Proof. Consider point $\underline{\bf y} = (y_2,...,y_d)$ with
\begin{equation}\label{n20}
| y_2 \cdots y_d| = \varepsilon,\,\,\,\,\,\,\, y_j >0,
\end{equation}We take  $(d-2)$ integers
$t_j, j = 2,...,d-1$  to satisfy
\begin{equation}\label{n1}
\left|
\log \frac{y_j}{\varepsilon^{1/(d-1)}} - t_j \log 2\right|\le \frac{\log 2}{2}
,
\,\,\,\,
2\le j \le d-1
\end{equation}
Then
\begin{equation}\label{n2}
\left|
\log \frac{y_d}{\varepsilon^{1/(d-1)}} - t_d \log 2\right|
=
\left|
\sum_{j=2}^{d-1}\left(\log \frac{y_j}{\varepsilon^{1/(d-1)}} - t_j\log 2\right)\right|
\le (d-1) \frac{\log 2}{2}
\end{equation}
 We see  from (\ref{n1},\ref{n2}) that
 $
 G_{-\underline{\bf t}} \, \underline{\bf y} \in \frak{B}_{\underline{\bf 0}} (\varepsilon) 
 $
 and so for any $\underline{\bf y}$ satisfying (\ref{n20})
 one has 
 $$
 \{\underline{\bf x}= (x_2,...x_d) \in \mathbb{R}^{d-1}:\,\,\,\,   |x_j|\le y_j,\,\,
 \,
 j = 2,...,d\} \subset  G_{\underline{\bf t}} \,  {B}_{\underline{\bf 0}} (\varepsilon)
 $$
 with $\underline{\bf t}$ defined in   (\ref{n1}).
 Also from (\ref{n1}) 
 and the triangle inequality  follows that if
$ \underline{\bf y} \in\frak{P}_{\varepsilon, X} $ then  integer $\underline{\bf t}$  defined in (\ref{n1}) satisfies
(\ref{maax}).$\Box$

 \vskip+0.3cm

Now  we consider function $\varphi (t)$ satisfying conditions 
  \begin{equation}\label{aaii}
  \varphi(t) \to 0\,\,\text{as}\,\, t\to\infty, \,\,\,\,\,\text{and}\,\,\,\,\,
  \varphi(s) \le 2\varphi (t)\,\,\text{ for}\,\,
  t\le s \le 2t,
  \end{equation}
  (every monotone $\varphi(t)$ satisfies (\ref{aaii}))
  and
  cover sets 
of the form
\begin{equation}\label{pab}
\frak{P}_{a,b}^i [\varphi] 
=
\left\{ {\bf x}\in \mathbb{R}^d: \,\,\,
x_i = \max_{j=1,...,d} |x_j| \in [2^a, 2^b],
\,\,\,
\Pi^d ({\bf x}) \le \varphi (x_i)
\right\}.
\end{equation}

 \vskip+0.3cm

In terms of function $ \varphi (\cdot )$ we define more specific parameters
$$
X(t_1) = 2^{t_1+1},\,\,\,\,\,\varepsilon(t_1) = \frac{2\varphi(2^{t_1})}{2^{t_1}},\,\,\,\,\, T(t_1) = T_{\varepsilon(t_1),X(t_1)},\,\,\,\,\,
\delta_j ({\bf t}) = 2^{(d-1)/2+t_j - t_1/(d-1)} (2\varphi(2^{t_1}))^{1/(d-1)}.
$$
From (\ref{suum}) and the choice of parameters  it follows  that 
 \begin{equation}\label{produkt1}
\prod_{j=2}^d \delta_j ({\bf t})  \ll \varepsilon(t_1) \ll \frac{\varphi(2^{t_1})}{2^{t_1}}
\end{equation}
and
 \begin{equation}\label{te1}
 T(t_1) \ll t_1 +\log \frac{1}{\varphi (2^{t_1})}
 .
\end{equation}
Consider boxes
$$
\frak{B}({\bf t}) =
\left\{ {\bf x} = (x_1,x_2,...,x_d) \in \mathbb{R}^d:\,\,\,
 2^{t_1} \le x_1 \le 2^{t_1+1},\,\,\, \underline{\bf x} \in 
\frak{B}_{\underline{\bf t}} (\varepsilon (t_1)).
\right\}
$$
We see that 
 \begin{equation}\label{produkt}
\frak{B}({\bf t}) =  I^{(1)} ({\bf t})\times  I^{(2)} ({\bf t})\times \cdots \times  I^{(d)} ({\bf t}),
\end{equation}
where 
$$
I^{(1)} ({\bf t}) = [2^{t_1}, 2^{t_1+1}],\,\,\,\,\,
I^{(j)} ({\bf t}) =[-\delta_j({\bf t}), +\delta_j({\bf t}),] ,\,\,\,  j = 2,...,d.
$$ 

\vskip+0.3cm
{\bf Lemma 4.}
{\it For domain (\ref{pab}) with integers $ a<b$  ($a$ large enough) and  $\varphi (\cdot)$ satisfying (\ref{aaii})
we have a covering by boxes
$$
\frak{P}^1_{a,b} [\varphi]
\subset
\bigcup_{t_1=a}^{b-1}\,\,\,\,\,
 \bigcup_{\underline{\bf t} }
  \frak{B} ({\bf t}),   
$$
where the inner  union is taken over all 
$\underline{\bf t}\in \mathbb{Z}^{d-1}$ such that  (\ref{suum})  holds and 
$
\max_{j=2,...,d-1} |t_j| \le T(t_1)$
}

\vskip+0.3cm
Proof
follows from Lemma 3, because for
$ 2^{t_1}\le x_1\le 2^{t_1+1}$
 the embedding
$$  \frak{P}_{\frac{\varphi(x_1)}{x_1}, x_1} \subset
\frak{P}_{\varepsilon (t_1), X(t_1)}
$$
is valid
due to  (\ref{aaii}):
for $ 2^{t_1}\le x_1 \le 2^{t_1+1}$ one has
$\frac{\varphi(x_1)}{x_1} \le 2\frac{\varphi (2^{t_1})}{2^{t_1}} = \varepsilon (t_1)$ and
$x_1\le X({t_1})$.
 So indeed,
$$
\frak{P}^1_{a,b} [\varphi]=
\bigcup_{t_1=a}^{b-1} \frak{P}^1_{t_1,t_1+1} [\varphi],
$$
and for fixed $t_1$ we have 
$$
\frak{P}^1_{t_1,t_1+1} [\varphi]=
\{ {\bf x}\in \mathbb{R}^d: 2^{t_1}\le x_1 \le 2^{t_1+1},\,\,\,
\underline{\bf x}\in   \frak{P}_{\frac{\varphi(x_1)}{x_1}, x_1} 
\}
\subset
\{ {\bf x}\in \mathbb{R}^d: 2^{t_1}\le x_1 \le 2^{t_1+1},\,\,\,
\underline{\bf x}\in   \frak{P}_{\varepsilon (t_1), X(t_1)}
\}\subset
$$
$$
\subset
\{ {\bf x}\in \mathbb{R}^d: 2^{t_1}\le x_1 \le 2^{t_1+1},\,\,\,
\underline{\bf x}\in    
 \bigcup_{\underline{\bf t} }
\frak{B}_{\underline{\bf t}} (\varepsilon (t_1))
\} =  \bigcup_{\underline{\bf t} }
  \frak{B} ({\bf t})\,\,\,\,\,\text{with}\,\,\,\,\,
  {\bf t} = (t_1,\underline{\bf t} ).
$$
Lemma 4 is proven.$\Box$


\section{$d$-dimensional lattices and the second application of Borel-Cantelli}

In addition to two fixed vectors $\pmb{e}_1, \pmb{e}_2$ defined in Section \ref{el} let us consider 
$d-2$ vectors 
\begin{equation}\label{eeee}
\pmb{e}_3=
\left(
\begin{array}{c}
e_{3,1}
\cr
e_{3,2}
\cr
\vdots
\cr
e_{3,d}
\end{array}
\right),
\,\dots, \,\,\,
\pmb{e}_d=
\left(
\begin{array}{c}
e_{d,1}
\cr
e_{d,2}
\cr
\vdots
\cr
e_{d,d}
\end{array}
\right)
\end{equation}
 in $\mathbb{R}^d$  related to coordinates $x_1,...,x_d$.
 Let 
$$
{E} = 
\left(
\begin{array}{ccc}
e_{3,1}&\cdots &e_{d,1}
\cr
\cdots&\cdots&\cdots
\cr
e_{3,d}&\cdots&e_{d,d}
\end{array}
\right) 
$$
be the corresponding  $d\times (d-2)$ matrix.
We consider matrix $E$ 
as an element of 
$d(d-2)$-dimensional space 
$
\mathcal{R} = \mathbb{R}^{d(d-2)} 
$.
For every matrix $E$ we define  a $d\times d$ matrix
$$
A({E}) = 
\left(
\begin{array}{ccccc}
e_{1,1}&e_{2,1}&e_{3,1}&\cdots &e_{d,1}
\cr
\cdots&\cdots&\cdots&\cdots&\cdots
\cr
e_{1,d}&e_{2,d}&e_{3,d}&\cdots&e_{d,d}
\end{array}
\right) 
$$
and the set 
\begin{equation}\label{teset}
\Lambda_{E} = \{ {\bf x} = z_1(\pmb{e}_1 +\beta \pmb{e}_2)+ z_2(\pmb{e}_2+\alpha \pmb{e}_1)+ 
 z_3\pmb{e}_3+\cdots+  z_d\pmb{e}_d,\,\,\,\ z_j \in \mathbb{Z},\,\,\, j = 1,...,d\},
\end{equation}
where $\alpha,\beta$ are the numbers fixed in Section \ref{s5s}.
Consider manifold
$$
\Sigma = \{ {E}\in \mathcal{R}:  {\rm det}\, A({E}) = 0\}.
$$
For every collection  
$$ 
u_{i,j} < v_{i,j} , i = 1,...,d; j = 3,..,d
$$
we define  $d(d-2)$-dimensional domains 
$$
\Omega^{(j)}
=
\{ (e_{3,j},...,e_{d,j}) \in \mathbb{R}^{d-2}:\,\,
u_{i,j} \le e_{i,j} \le v_{i,j},  i = 3,...,d\},\,\,\,\, j = 1,...,d  
$$
and
$$
\Omega
=
\Omega^{(1)}\times \cdots \times \Omega^{(d)}
\subset \mathcal{R}.
$$
Under the assumption
\begin{equation}\label{kompas}
\Omega \cap \Sigma = \varnothing,
\end{equation}
for all $E\in \Omega$ vectors $ \pmb{e}_1, \pmb{e}_2,  \pmb{e}_3,...,  \pmb{e}_d$
are linearly independent.  
So in this case $\Lambda_E$ from (\ref{teset}) is a $d$-dimensional lattice.
  In such a case
every vector ${\bf x} \in \Lambda_E$ may be represented as
\begin{equation}\label{bisa}
{\bf x} =
{\bf x}({\bf z})
= (x_1,...,x_d) =  z_1(\pmb{e}_1+\beta \pmb{e}_2)+
z_2(\pmb{e}_2+\alpha\pmb{e}_1) + z_3\pmb{e}_3
+\cdots+  z_d\pmb{e}_d,\,\,\,\,  \text{with}
\,\,\,\,
{\bf z} = (z_1,...,z_d)\in \mathbb{Z}^d.
\end{equation}
A point of $\Lambda_E$ belongs to $\frak{L}$ if and only if $z_3=...=z_d=0$, so 
$\Lambda_E\cap\frak{L} = \Gamma_{\frak{L}}$.
For ${\bf x}$ 
(or ${\bf z}$) from the equality above we consider values 
\begin{equation}\label{define}
X = \max_{j = 1,..,d} |x_j|\,\,\,\,\,\,\,\,
\text{and}
\,\,\,\,\,\,\,\,
Z = \max_{j = 1,..,d} |z_j|.
\end{equation}
In addition, for   ${\bf z} = (z_1,...,z_d)\in \mathbb{Z}^d$ we consider short vectors 
$${\bf z}_- = (z_1,z_2)\in \mathbb{Z}^{2}
\,\,\,\,\,\,\,\,\, \text{and}
\,\,\,\,\,\,\,\,\,
{\bf z}_+ = (z_3,...,z_d)\in \mathbb{Z}^{d-2}.
$$

\vskip+0.3cm

{\bf Second Main Lemma.} {\it Let $\Omega$ satisfy (\ref{kompas}).
Then there exists $K= K(\Omega)\ge 1$ with the following property.
  Assume that  function $\varphi (t)$  
  satisfies (\ref{aaii})
  and both  series 
\begin{equation}\label{sse}
\sum_{\nu=1}^\infty \nu^{d-1}\cdot\varphi\left(\frac{2^\nu}{K}\right),\,\,\,\,\,\,\,\,\,\,\,
\sum_{\nu=1}^\infty \, \nu\cdot
\left(\log \frac{1}{\varphi (K 2^{\nu})}\right)^{d-2}\!\!\cdot\varphi\left(\frac{2^\nu}{K}\right)
\end{equation}
 converge. 
Then
for almost all 
$
{E}  \in 
\Omega   
$ 
  there exists 
  $\eta  = \eta (E)>0$ such that for all  
$
{\bf x} \in \Lambda_{E}
$
with 
\begin{equation}\label{zedplus}
Z_+ = |{\bf z}_+ |=
\max_{j= 3,...,d} |z_j| \neq 0
\end{equation}
we have
\begin{equation}\label{f1}
\Pi^d({\bf x}) \ge  \eta \varphi (X).
\end{equation}
}

\vskip+0.3cm
 Proof.  
 We consider the  "exceptional" set 
$$
\mathcal{M} = 
\mathcal{M}  [\varphi] =
\{
E\in \mathcal{R}:
\text{there are infinitely many }\,\,\,
{\bf x} \in \Lambda_E \cap \frak{P} [\varphi] \,\,\,
\text{with}\,\,\, Z_+ \neq 0 
\}
$$
and prove that  for any $\Omega$ satisfying (\ref{kompas})
there exists $K(\Omega)$  which depends on $\Omega$ but does not depend on $\varphi(\cdot)$, such that  under the convergence condition
we get
$
\mu (\mathcal{M}\cap \Omega) =0.
$
This will prove Second Main Lemma.

Define
$$
\mathcal{W}_\nu =\{ E\in \Omega:
{\bf x}({\bf z})\in \frak{P} [\varphi] \,\,\,\,
\text{for some}\,\,\,\, {\bf z}\,\,\,\,
\text{with}
 \,\,\, Z_+ \neq 0 \,\,\, \text{and}
\,\,\,\,
 2^\nu \le Z\le 2^{\nu+1}\}.
$$
As
$$
\mathcal{M}  = \bigcap_{\nu_0=1}^\infty \bigcup_{\nu = \nu_0}^\infty \mathcal{W}_\nu
$$
by Borel-Cantelli argument it is enough to prove that the series
\begin{equation}\label{meas}
\sum_{\nu=1}^\infty \mu (\mathcal{W}_\nu)
\end{equation}
converges.
So we 
calculate an upper bound for the measure $\mu (\mathcal{W}_\nu)$.

\vskip+0.3cm

 From compactness of $\Omega$ and 
  (\ref{kompas}) we see that   for $X$ and $Z$ defined in (\ref{define}) one has 
$ X\asymp_\Omega Z$. By the same reason, 
there exist   a positive integer 
$h =  h(\Omega)$ such that  for any $\nu$ and for any
${\bf z}$ with 
\begin{equation}\label{zeta}
 2^\nu \le Z\le 2^{\nu+1}
 \end{equation}
  for the corresponding value of $X$ 
inequalities 
\begin{equation}\label{kompas1}
 2^{\nu-h}\le X\le 2^{\nu+h}
\end{equation}
are  valid
for any choice of vectors (\ref{eeee}) with  $ E \in \Omega$.
In particular for ${\bf z}$  satisfying (\ref{zeta}) for any $ E \in \Omega$  we conclude that 
\begin{equation}\label{kompas12}
{\bf x}({\bf z})\in \frak{P} [\varphi] \,\,\,\text{with }\, {\bf z}\, \text{satisfying (\ref{zeta})}
\,\,\,\,
\Longrightarrow
\,\,\,\,
{\bf x}({\bf z})\in
\bigcup_{i=1}^d
\frak{P}_{\nu - h,\nu+h}^i [\varphi] \,\,\,\,\,\text{or}\,\,\,\,\,
{\bf x}(-{\bf z})\in
\bigcup_{i=1}^d
\frak{P}_{\nu - h,\nu+h}^i [\varphi] 
 .
\end{equation}
Define
$$
\mathcal{W}_\nu' =\left\{ E\in \Omega:
{\bf x}({\bf z})\in
\bigcup_{i=1}^d
\frak{P}_{\nu - h,\nu+h}^i [\varphi] \,\,\,\,
\text{for some}\,\,\,\, {\bf z}\,\,\,\,
\text{satisfying (\ref{zeta}) with}
 \,\,\, Z_+ \neq 0 \right\},
$$
$$
\mathcal{W}_\nu'' =\left\{ E\in \Omega:
{\bf x}({\bf z})\in
\frak{P}_{\nu - h,\nu+h}^1 [\varphi] \,\,\,\,
\text{for some}\,\,\,\, {\bf z}\,\,\,\,
\text{satisfying (\ref{zeta}) with}
 \,\,\, Z_+ \neq 0 \right\}
$$
and
$$
\mathcal{V}_{\nu, {\bf t} }
=\left\{ E\in \Omega:
{\bf x}({\bf z})\in
 \frak{B} ({\bf t})\,\,\,\,
\text{  for some}\,\,\,\, {\bf z}\,\,\,\,
\text{satisfying (\ref{zeta}) with}
 \,\,\, Z_+ \neq 0 \right\}
 $$
 (without loss of generality  we can hold estimates for $i=1$).
 The set $\frak{B} ({\bf t})$ can be decomposed as in (\ref{produkt}). So if ${\bf x}={\bf x}({\bf z}) \in \frak{B} ({\bf t})$ is of the form (\ref{bisa}), we conclude that 
 $$
 x_j \in  I^{(j)} ({\bf t}) ,\,\,\,\,\,\,\, j = 2,....,d,
 $$
 or for fixed ${\bf z}, {\bf t}$ and any $j =2,...,d$  the values of $e_{i,j}$ satisfy
 $$
 |z_1 (e_{1,j}+\beta e_{2,j})+ z_2(e_{2,j}+\alpha e_{1,j})+ z_3e_{3,j}+...+z_d e_{d,j} |\le \delta_j ({\bf t}).
 $$
 So, being ${\bf z}$ and ${\bf t}$  fixed, condition ${\bf x}={\bf x} \in \frak{B} ({\bf t})$ leads to 
 $$
E\in \Omega_{{\bf z},{\bf t}}= \Omega^{(1)}\times \Omega_{{\bf z},{\bf t}}^{(2)}\times \cdots \times \Omega_{{\bf z},{\bf t}}^{(d)}
\subset \mathcal{R},
$$
where 
\begin{equation}\label{aom}
 \Omega_{{\bf z},{\bf t}}^{(j)}  =\Omega^{(j)}  \cap \{ (e_{3,j},...,e_{d,j}) \in \mathbb{R}^{d-2}:
|z_1 (e_{1,j}+\beta e_{2,j})+ z_2(e_{2,j}+\alpha e_{1,j})+z_3e_{3,j}+...+z_d e_{d,j} |\le \delta_j ({\bf t})\},\,
j=2,...,d.
\end{equation}
Finally, we get the inclusion
\begin{equation}\label{inka}
\mathcal{V}_{\nu, {\bf t} }
\subset \bigcup_{{\bf z}:\,\, 2^\nu \le Z\le 2^{\nu+1},\,  Z_+\neq 0}
\Omega_{{\bf z},{\bf t}}.
\end{equation}
Implication (\ref{kompas12})  together with definitions of 
$\mathcal{W}_\nu',
\mathcal{W}_\nu'',
\mathcal{V}_{\nu,{\bf t}}$,
 Lemma 4
 and (\ref{inka})
 give us inequalities
$$
\mu(\mathcal{W}_\nu) \le 
\mu(\mathcal{W}_\nu') \le
d \cdot 
\mu(\mathcal{W}_\nu '')
 \ll  
  \sum_{t_1=\nu - h,}^{\nu + h,}\,\,\,\,\,
 \sum_{\underline{\bf t} }\,\,
 \mu ( \mathcal{V}_{\nu, {\bf t} })\le
  \sum_{t_1=\nu - h,}^{\nu + h,}\,\,\,\,\,
 \sum_{\underline{\bf t} }\,\,\,\,\,\,\,\,
 \sum_{{\bf z}:\,\, 2^\nu \le Z\le 2^{\nu+1},\,  Z_+\neq 0}\,\,\,
 \mu (
\Omega_{{\bf z},{\bf t}}
 )
 ,
$$
where the  sum 
$\sum_{\bf t}$
is taken over all 
$\underline{\bf t}\in \mathbb{Z}^{d-1}$ such that  (\ref{suum})  holds and 
$
\max_{j=2,...,d-1} |t_j| \le T(t_1)$.

We note that by (\ref{aaii}), for
$\nu-h \le t_1 \le \nu+h$  we have 
$4^{-h} \varphi(2^{\nu+h} \le \varphi (2^{t_1}) \le 4^h \varphi (2^{\nu-h})$
 and hence
 $T(t_1) \le T(\nu+h)+2h, \varepsilon (t_1) \le 4^h \varepsilon (\nu-h)$.
Taking into account upper bound (\ref{te1}) we get
 \begin{equation}\label{aaa}
 \mu(\mathcal{W}_\nu) \ll_\Omega
 (T(\nu+h))^{d-2} 
  \sum_{{\bf z}} \max_{\bf t}
   \mu (
\Omega_{{\bf z},{\bf t}}
 )
\ll_\Omega
 \left(
 \nu + \log \frac{1}{\varphi (2^{\nu+h})}
 \right)^{d-2}
 \sum_{{\bf z}:\,\, 2^\nu \le Z\le 2^{\nu+1},\,  Z_+\neq 0}\,\,\,
 \max_{\bf t}\, 
 \mu (
\Omega_{{\bf z},{\bf t}}
 ),
 \end{equation}
 where $\max_{\bf t}$ 
is taken over all 
${\bf t} = (t_1,\underline{\bf t}) \in \mathbb{Z}^{d}$ such that
$ \nu-h\le t_1 
\le \nu+h$,
  (\ref{suum})  holds and 
\begin{equation}\label{aaaa}
\max_{j=2,...,d-1} |t_j| \le  K_1(\Omega)  \left(
 \nu + \log \frac{1}{\varphi (2^{\nu+h})}
 \right)
 \end{equation}
with positive constant $K_1(\Omega)  $ depending on $\Omega$. However, the bound (\ref{aaaa}) will be not of importance for our final consideration.

\vskip+0.3cm
We need two more easy lemmas. The proofs will be given  in Section \ref{sle}.

\vskip+0.3cm
 {\bf Lemma 5.}
 {\it 
 For fixed vector ${\bf z}_+$ the number of two-dimensional vectors ${\bf z}_-=(z_1,z_2)$ with 
 $$ Z_-=\max (|z_1|,|z_2|) \le Q\,\,\,\,\,\text{ and }\,\,\,\,\,
 \Omega_{{\bf z},\bf{t}}^{(j)}  \neq \varnothing
 {\,\,\,  \forall \, j = 2,...,d}
 $$
   is $\ll_\Omega QZ_+$ where $Z_+$ is defined in (\ref{zedplus}).}
 \vskip+0.3cm

\vskip+0.3cm
 {\bf Lemma 6.}
 {\it   Assume that $Z_+\neq 0$. Then for any $ j =2,...,d$  for $(d-2)$-dimensional Lebesgue measure 
 $\Omega_{{\bf z},\bf{t}}^{(j)} \subset \mathbb{R}^{d-2}$ the 
upper bound 
 $$
 \mu_{d-2} ( \Omega_{{\bf z},\bf{t}}^{(j)} ) \ll_\Omega \frac{ \delta_j ({\bf t})}{Z_+}
 $$
 holds uniformly in $\bf t$.
 
 }
\vskip+0.3cm


We continue estimating $  \mu(\mathcal{W}_\nu) $. 
The argument below, as well as the proof of 
  Lemmas  5 and 6 were utilised recently by the author in \cite{Uri},  Claims 7.6.1, 7.6.2.

Lemma 6 together with (\ref{produkt1}) gives the  bound
$$
\max_{\bf t}
 \mu (
\Omega_{{\bf z},{\bf t}} )  \ll_\Omega   \max_{\bf t} \frac{ \prod_{j=2}^d\delta_j ({\bf t})}{Z_+^{d-1}}\ll
\frac{\varepsilon(\nu-h)}{Z_+^{d-1}}
 \ll_\Omega  \frac{\varphi\left(2^{\nu-h}\right)}{2^\nu Z_+^{d-1}}.
$$
Now
$$
 \sum_{{\bf z}:\,\, 2^\nu \le Z\le 2^{\nu+1},\,  Z_+\neq 0}\,\,\,
 \max_{\bf t}\, 
 \mu (
\Omega_{{\bf z},{\bf t}}
 )\ll_\Omega
 \frac{\varphi\left(2^{\nu-h}\right)}{2^\nu }.
\sum_{{\bf z}_+: 0< Z_+\le 2^{\nu+1}}
\frac{1}{Z_+^{d-1}}
\sum_{
\begin{array}{c}
{\bf z}_-: Z_- \le2^{\nu+1} \,\,\,\text{and}\cr
\Omega_{{\bf z},\bf{t}}^{(j)}  \neq \varnothing
 \,\,\,  \forall \, j = 2,...,d
\end{array}
} 1\ll
$$
$$
\varphi (2^{\nu-h}) \cdot
\sum_{{\bf z}_+: 0< Z_+\le 2^{\nu+1}}
\frac{1}{Z_+^{d-2}}
= \varphi (2^{\nu-h}) \cdot
\sum_{q=1}^{2^{\nu+1}}  \frac{\#\{ {\bf z}_+\in \mathbb{Z}^{d-2}: Z_+ = q\} }{q^{d-2}}\ll 
\nu \cdot \varphi (2^{\nu-h})
,
$$
  because of Lemma 5  with $Q = 2^{\nu+1}$ and inequality $ \#\{ {\bf z}_+\in  \mathbb{Z}^{d-2}: Z_+ = q\} \asymp q^{d-3}$. Finally,  continuing the estimate from (\ref{aaa}) we get 
  $$
 \mu(\mathcal{W}_\nu) \ll_\Omega
\nu \cdot \left( \nu^{d-2}+\left(\log \frac{1}{\varphi (2^{\nu+h})}\right)^{d-2}\right)\cdot \varphi (2^{\nu-h}).
$$
Put $K(\Omega) = 2^{h}$. Then  the assumption that series (\ref{sse}) converge shows that 
the series from (\ref{meas}) converges  also, and Second Main Lemma is proven.$\Box$

 \section{Exact approximations and proof of Theorem 2}

 Define 
 $$
 \tilde{w} (t)= \left(\frac{16
 }{15}\right)^{d-1} \sigma^{-1} \, \frac{w\left( \frac{15}{16}  \max_{2\le j \le d}|e_{2,j}| \cdot t\right)}{ t^{d-1}},
 $$ 
 with $\sigma = \sigma(\frak{L})$ defined in (\ref{ito1}).
 For $\delta >0$ we take irrational $\alpha$ satisfying the following properties:
 
  \vskip+0.3cm
  \noindent
 (A) $\alpha$ satisfies (\ref{nee}) and $\left|\alpha - \frac{1}{4}\right|\le\delta $;
 
   \vskip+0.3cm
   \noindent
 (B) there exists an infinite  sequence of convergents $\frac{p_{\nu_k}}{q_{\nu_k}}, k=1,2,3,...$ to $\alpha$ such that 
 $
 ||q_{\nu_k}\alpha|| \le  { \tilde{w}\left(   q_{\nu_k}\right)}
  $;
 
  \vskip+0.3cm
   \noindent
 (C) for every $\nu$ large enough  
 $
 ||q_{\nu}\alpha|| \ge (1-\delta)   { \tilde{w}\left(  q_{\nu}\right)}.
 $

 \vskip+0.3cm
 
 Existence of such $\alpha$ follows for example from Theorem B mentioned in Section 2, because of $ t\cdot  \tilde{w} (t) \downarrow 0$ as $ t \to \infty$.
 
 Then we apply Lemma 2 with $\psi(t) = t^{-2}$ and take $\beta \in \mathcal{B}$ with 
 \begin{equation}\label{beeta}
 \left|\beta - \frac{1}{4}\right|\le \delta.
  \end{equation}
  Now we 
  take $\Omega$ to satisfy (\ref{kompas}) and consider constant $K= K(\Omega)$ from the Second Main Lemma. 
  
  We need one more easy lemma. A proof will be given in the next section.
  
  \vskip+0.3cm
 {\bf Lemma 7.} {\it 
 Let $w(t)$ satisfy conditions of Theorem 2 and $K\ge 1$ be fixed. Then there exist function $\varphi(t)$ 
 satisfying (\ref{aaii})
  and all the  conditions of the Second Main Lemma  with the given value of $K$ and in addition
\begin{equation}\label{finalist}
  \lim_{t\to \infty} \frac{\varphi(t)}{w(t)} = +\infty .
  \end{equation}}
    \vskip+0.3cm
    
   We apply the Second Main Lemma  for function $\varphi (t)$ from Lemma 7 and take arbitrary $E$ satisfying its conclusion. 
Then  lattice $\Lambda_E$ satisfies the conclusion of Theorem 2. Indeed, consider  several cases.

 \vskip+0.3cm
 \noindent
  1) If ${\bf x}\in \Lambda_E\setminus\frak{L}$ by the Second Main Lemma we get inequality 
  (\ref{f1}) which  together with (\ref{finalist}) gives 
  $\Pi^d({\bf x} ) \ge \eta \varphi (X) \ge  w(X)$ for all $X$ large enough.
  
 \vskip+0.3cm

   \noindent
  2) If   ${\bf x}\in \Gamma_{\frak{L}}\subset \Lambda_E $ but not of the form 
  ${\bf z}_{m,n} $ with $ (m,n) = s(q_\nu, -p_\nu)$, by the First Main Lemma we get 
   $\Pi^d({\bf x} )\gg_{\alpha,\beta}  1 $ and so again  $\Pi^d({\bf x} )\ge w(X)$ for all $X$ large enough.
  
   \vskip+0.3cm
   \noindent
    3) If   ${\bf x}= {\bf z}_{m,n} = s (q_\nu, -p_\nu)$,  $ \nu = \nu_k$ then by the First Main Lemma 
    conditions (\ref{beeta}) and (A), (B), (C)  
    for $s= \pm 1$ we have 
     $$
  \Pi^d ({\bf x})  
  = \sigma \left(\frac{15}{16}\right)^{d-1}q_{\nu+k}^{d-1}||q_{\nu_k}\alpha || \cdot (1+O(\delta))
  =\sigma \left(\frac{15}{16}\right)^{d-1}q_{\nu_k}^{d-1}\tilde{w}\left(  q_{\nu_k}\right) \cdot (1+O(\delta))=
  w\left( \frac{15}{16}  \max_{2\le j \le d}|e_{2,j}| \cdot  q_{\nu_k}\right)  \cdot (1+O(\delta)).
 $$
 Taking into account (\ref{1x})  and condition (\ref{www00}) we continue with 
    $$
  \Pi^d ({\bf x})  
  = w\left(
  X \cdot (1+O(\delta)
  \right)  \cdot (1+O(\delta))= w(X) \cdot (1+O(\varepsilon)),
  $$
for $\delta = \delta (\varepsilon)$ small enough and $\nu_k$ large enough.

       \vskip+0.3cm
   \noindent
      4) If   ${\bf x}= {\bf z}_{m,n} = s(q_\nu, -p_\nu)$,  $ \nu = \nu_k$  but $ s\neq \pm 1,0$, analogously to the case 3) we have 
          $$
  \Pi^d ({\bf x})  
  \ge w\left(
  X
  \right)  \cdot (2+O(\varepsilon)) ,
  $$
   for $\delta = \delta (\varepsilon)$ small enough and for $\nu_k$ large enough.  
      
      \vskip+0.3cm
   \noindent
      5) If   ${\bf x}= {\bf z}_{m,n} = s(q_\nu, -p_\nu)$ for $ \nu \neq \nu_k$ and $ s\neq 0$,
      by condition (C)   analogously to the case 3) we have 
    $$
  \Pi^d ({\bf x})  
  \ge w\left(
  X
  \right)  \cdot (1-O(\varepsilon)) ,
  $$
   for $\delta = \delta (\varepsilon)$ small enough and for all $\nu$ large enough.

\vskip+0.3cm
Theorem 2 is proven.$\Box$

\section{ Proofs of Lemmas 5, 6, and 7}\label{sle}

To prove Lemma 5 
denote
$
\lambda = z_3e_{3,j}+...+z_d e_{d,j} .
$
Consider the value  $e_{2,j}+\alpha e_{1,j}$. Since $\frak{l}_j =\{ x_j =0\}$, we have 
$\xi_j = -e_{1,j}/e_{2,j}$, so $e_{2,j}+\alpha e_{1,j} =e_{2,j}(1-\alpha \xi_j) \neq 0$
 by(\ref{eee}) and (\ref{nee}).
For any fixed value of $z_1$  condition $ \Omega_{{\bf z},\bf{t}}^{(j)}  \neq \varnothing$ leads to
$$
  \left| z_2 + \frac{z_1(e_{1,j}+\beta e_{2,j})}{e_{2,j}+\alpha e_{1,j}}\right| \le \frac{ |\lambda|+ \delta_j ({\bf t})}{|e_{2,j}+\alpha e_{1,j}|}.
$$
We can choose $j $ such that $ \delta_j ({\bf t}) \ll 1$.
For any choice of $E\subset \Omega$  we see that 
$|\lambda |\ll_\Omega Z_+$.
 So for any fixed value of $z_1$  we get
 not more than $ O_\Omega (Z_+)$ possible values for integer $z_2$.
At the same time $z_1$ admit $O(Q)$ possible values.
Lemma 5 is proven.$\Box$

 \vskip+0.3cm

Lemma 6 immediately follows from the observation that 
 $ \Omega_{{\bf z},\bf{t}}^{(j)} $  is a strip in 
 the box $ \Omega^{(j)}\subset \mathbb{R}^{d-2}$ of width $ 
 \frac{2 \delta_j ({\bf t})}{ \sqrt{z_{3}^2+...+z_d^2}} \ll  \frac{2 \delta_j ({\bf t})}{Z_+}$.
Lemma 6 is proven.$\Box$

\vskip+0.3cm


\vskip+0.3cm

To prove Lemma 7 one should observe that from the assumptions of Theorem 2 we conclude that the series
$$
\sum_{\nu=1}^\infty \nu^{d-1} w\left(\frac{2^\nu}{K}\right)
$$
converges.
Then one can define  a strictly increasing sequence of positive integers   $\nu_l, l=1,2,3,...$ such that 
$$
\sum_{\nu=\nu_l}^\infty \nu^{d-1} w\left(\frac{2^\nu}{K}\right)\le \frac{1}{4^l}.
$$
We put
$$
\lambda(t) = 2^l\,\,\,\,\, \text{for}\,\,\,\,\,
{2^{\nu_l}} \le t < {2^{\nu_{l+1}}};\,\,\,\,\,\,
\text{so} \,\,\,\,\,\, \lambda(t) \to \infty, t\to \infty
,
$$
and define
$$
\varphi(t) = w(t) \lambda(Kt).
$$
Then (\ref{aaii}) and (\ref{finalist})  are satisfied and 
$$
\sum_{\nu=\nu_l}^{\nu_{l+1}-1} \nu^{d-1} \varphi \left(\frac{2^\nu}{K}\right)\le
\sum_{\nu=\nu_l}^{\nu_{l+1}-1} \nu^{d-1} w \left(\frac{2^\nu}{K}\right) \lambda (2^\nu)\le
\frac{1}{2^l},
$$
so the first series from (\ref{sse}) converges.

Then putting $a_\nu =\varphi \left(\frac{2^\nu}{K}\right) , r= d-1$
 we observe that   
$$
\nu a_\nu^{3/4} = (\nu^r a_\nu)^{3/4} \cdot \nu^{1- 3r/4}\,\,\,\,\,\,
\text{and}
\,\,\,\,\,\,
(\nu^{1- 3r/4})^4 = \nu^{4-3r} \le \nu^{-2}.
$$
From Hölder inequality 
$$
\sum_\nu u_\nu v_\nu \le \left(\sum_\nu u_\nu^{4/3} \right)^{3/4}
\left(\sum_\nu v_\nu^{4} \right)^{1/4}
 $$
 and convergence of the first series from (\ref{sse})  we see that $\sum_\nu \nu a_\nu^{3/4}$ converges.
 We use  inequality 
 (\ref{www0000}) which is  a corollary of (\ref{www00})  to see that the second series from (\ref{sse}) converges also.
 Lemma 7 is proven.$\Box$

 \vskip+0.3cm
 
 {\bf AI declaration.} The author used Claude Anthropic (model Fable 5.1) to determine inaccuracies in the final version of the manuscript.

\end{document}